\documentclass[12pt]{article}
\usepackage{amssymb,amsmath}
\begin{document}
\title{Branched immersions and braids}
\author{Marina Ville}
\date { }
\maketitle
\newtheorem{slem}{Sublemma}
\newtheorem{lem}{Lemma}
\newtheorem{thm}{Theorem}
\newtheorem{prop}{Proposition}
\newtheorem{defn}{Definition}
\newtheorem{cor}{Corollary}

ABSTRACT. Branch points $p$ of a real $2$-surface $\Sigma$ in a
$4$-manifold $M$  generalize branch points of complex curves in
complex surfaces: for example, they can occur as singularities of
minimal surfaces. We investigate such a branch point $p$ when
$\Sigma$ is topologically embedded. It defines a link $L(p)$, the
components of which are closed braids with the same axis up to
orientation. If $\Sigma$ is closed without boundary, the
contribution of $p$ to the degree of the normal bundle of $\Sigma$
in $M$ can be computed on the link $L(p)$, in terms of the algebraic
crossing numbers of its components and of their linking numbers with
one another.


\section{Introduction}
\subsection{Motivation}
We investigate here the generalization to the real
case of a construction which is well understood in the complex algebraic case.
Consider a complex algebraic curve $C$ in $\Bbb{C}^2$ which possesses a branch
point $p$. For a small enough positive real number $\epsilon$ the intersection
of $C$ with the sphere $\Bbb{S}_\epsilon(p)$ centered at $p$ and of radius $\epsilon$
is a link which is called an {\it algebraic} link. The link type does not depend
on $\epsilon$. Numerical invariants  of this link (e.g. its genus if it is a
knot) can be read on the singularity of $C$ at $p$ and vice versa.\\

Complex branch points generalize to real branch points; in
particular these are the singularites of minimal surfaces in
Riemannian manifolds. They are non generic singularities of real
surfaces which lowest order term is similar to the lowest
order term of a complex branch point.\\

In the present paper we consider a closed surface $S$ which is immersed in a $4$-manifold $M$ except at a finite number of branch points.\\
For such a surface we can define a tangent and a normal bundle, which we denote $TS$ and $NS$. We will give a
precise definition below; let us just say for now that $TS$
is about the intrinsic topology of the surface and $NS$ reflects
its extrinsic topology (how it sits inside $M$).\\

We make an extra assumption which is always satisfied in the complex
case, but not in the real case, namely that the singularities of $S$
are isolated (even minimal surfaces can have
real codimension $1$ singularities).\\

If $p$ is a singular point of $S$, we consider the sphere
$\Bbb{S}_\epsilon(p)$ centered at $p$ and of radius $\epsilon$; the
intersection of $S$ and $\Bbb{S}_\epsilon(p)$ is
a link which we denote $L^\epsilon$. \\

If there is a single branched disk going through $p$, $L^\epsilon$
is a knot which is transverse to the standard contact structure on
$\Bbb{S}_\epsilon(p)$; moreover it is immediately available in the
form of a {\it closed braid} as defined by Bennequin. If there are
several different disks going through $p$, there is no obvious
contact structure to which $L^\epsilon$ is transverse. However we
can find an axis w.r.t. which and up to orientation the components
of $L^\epsilon$ are closed braids. The braid indices of the
components and their algebraic crossing numbers do not depend on
$\epsilon$;
and neither do the linking numbers of one component with another. \\

These quantities have an interpretation in terms of the topology of
$S$. We can derive from  the braid index of $L^\epsilon$ the
contribution of $p$ to the degree of the tangent bundle $TS$. This
has been known for a long time. In the present paper we show that
the algebraic crossing numbers and linking numbers of the components
of $L^\epsilon$ give the
contribution of $p$ to the degree of the normal bundle $NS$.\\

As an application we can recover the (already known) formula for the
self-linking number of iterated torus knots. \\
\\

The starting point of this work was a conversation with Alexander
Reznikov many years ago. He gave me advice, friendship and
encouragement. Not long ago he died tragically. This paper is
inscribed to his memory.
\\
\\
ACKNOWLEDGEMENTS. The author would like to thank Jim Eells for helpful advice
and conversation, Denis Auroux for providing her with much needed information
about $4$-dimensional topology and Marc Soret for introducing her to braids.
\subsection{Preliminaries}
\subsubsection{Transverse knots - Closed braids}
We consider $\Bbb{S}^3$ as the unit sphere in the complex plane $\Bbb{C}^2$.
The
complex structure $J$ on $\Bbb{C}^2$ enables us to define a {\it contact
structure}
$\xi$ on $\Bbb{S}^3$. Namely, if  $q$ is a point in $\Bbb{S}^3$, the plane
$\xi(q)$ is the plane
tangent to
$\Bbb{S}^3$ at $q$ and orthogonal to $Jq$. Stereographic projection maps this
to the standard contact structure in $\Bbb{R}^3$ which we also denote $\xi$:
in other words, up to isomorphism,
$\xi$ is the contact structure for which the contact planes are the kernels
of the $1$-form written in
cylindrical
coordinates as
$$\rho^2d\theta+dz.$$

A {\it transverse knot} (in $\Bbb{S}^3$ or $\Bbb{R}^3$) is a smooth
map $\gamma$ from the circle $\Bbb{S}^1$ to $\Bbb{S}^3$ or
$\Bbb{R}^3$ such that
$$\forall t\in\Bbb{S}^1,\ \ \ \ \gamma'(t)\notin\xi(\gamma(t))$$
i.e. the tangent vectors to the knot never belong to the contact planes. In
view of the $1$-form given above, this means that if the knot is written in
cylindrical coordinates $(\theta,\rho,z)$ it verifies
$$\forall t\in \Bbb{S}^1, \frac{z'(t)}{\theta'(t)}\neq-(\rho(t))^2.$$

Two transverse knots are {\it transversally isotopic} if they have the same knot type and moreover they are isotopic through transverse knots.\\
A special case of transverse knots (in $\Bbb{R}^3$) is given by what Bennenquin
calls {\it closed braids}:
\begin{defn}
A knot $K$ in $\Bbb{R}^3$ is a {\bf closed braid} if, when written in
cylindrical
coordinates
$$t\mapsto (\rho(t),\theta(t), z(t))$$
it verifies
$$\forall t\in\Bbb{S}^1, \rho(t)\neq 0, \theta'(t)>0.$$
\end{defn}
In this definition the $z$-axis is called the {\it axis} of the braid.
\\Bennequin ([Be 1]) proved that every transverse knot is transversally isotopic
to a
closed braid.\\

There is a similar definition for closed braids in the $3$-sphere;
in this case the axis of the braid is an oriented great circle. For
further use we state
\begin{lem}
Let $L$ be a link in $\Bbb{S}^3$. Let $P$ be a plane in $\Bbb{R}^4$ and let
$(\epsilon_1, \epsilon_2)$ be a (not necessarily orthonormal) basis of the
orthogonal complement $P^\perp$. \\

We denote the orthogonal projection of $L$ to $P^\perp$ by
$$x_1(t)\epsilon_1+x_2(t)\epsilon_2.$$
The following two assertions are equivalent\\
1) up to orientaton the great circle $\Gamma$ in $P$ is a braid axis for $L$\\
2) the projection of $L$ to $P^\perp$ verifies
$$x_1^2(t)+x_2^2(t)\neq 0,\ \ \  x_1(t)x'_2(t)-x_2(t)x'_1(t)\neq 0.$$
\end{lem}
We conclude this section by recalling a few invariants for a  closed braid $K$. For more details see for example [B-W].
The {\it braid index} $n(K)$ of $K$ is the linking number of $K$ with the
oriented $z$-axis. A generic projection of $K$ onto the $(\rho,\theta)$ plane
has only transverse double crossing  points; it is called  a {\it closed braid
projection}. We
assign a sign to such an intersection point in the following manner:
consider a basis $(u_1,u_2)$ of $\Bbb{R}^2$ where the $u_1$ (resp. $u_2$) is the vector tangent to the strand of $K$ which is on top (resp.
on the bottom).
If $(u_1,u_2)$ basis is a positive (resp. negative) basis the point will be counted positively (resp. negatively).\\

The {\it algebraic crossing number} $e(K)$ of $K$ is the signed
number
crossing points of a closed braid projection of $K$.\\

The quantity
$$sl(K)=n(K)-e(K)$$
is called the {\it self-linking number} of the transverse knot. Bennequin
introduced it in [Be 1] and proved that it is an invariant of transverse
isotopy.

\subsubsection{Branched immersions}
We recall a few (standard) definitions. For more details on branched immersions
we refer the reader to [G-O-R].
\begin{defn}
A {\bf branched disc} in a smooth manifold $M$ is a map from a disc $D$ centered at $0$ in $\Bbb{C}$ to $M$
which is an immersion except at $0$ and which writes in a neighbourhood of $0$\\
$f^1(z)=Re(z^N) + o_1(|z|^N)$\\
$f^2(z)=Im(z^N)+o_1(|z|^N)$\\
$f^k(z)=o_1(|z|^N)$ for $k>2$\\
where $z$ is a local isothermic coordinate on $D$ around $0$, the $f^i(z)$'s
are the coordinates of $f(z)$ in some well-chosen chart on $M$ around $f(0)$ and $N$ is an
integer, $N\geq 2$.\\

The notation $o_1(|z|^N)$ means that the function is an $o(|z|^N)$
and its
first derivatives are $o(|z|^{N-1})$'s.\\

The integer $m=N-1$ is called the {\bf branching order} of $f$ at
$0$.
\end{defn}
\begin{defn}
A map $f:\Sigma\longrightarrow M$ from a Riemann surface $\Sigma$ to a
smooth manifold $M$ is a {\bf branched immersion} if it is an immersion except at
a discrete set of points called {\bf branch points} in a neighbourhood of which $f$ is parametrized by branch discs.
\end{defn}
Note that there can be several branched discs going through the same branch
point.\\

If $f$ is a branched immersion as in Def. 3  one can check that the
map from $\Sigma$ to the Grassmannian of oriented $2$-planes
$G_2^+(M)$ given by
$$x\mapsto f_*(T_x\Sigma)$$
(where $T_x\Sigma$ is the tangent plane to $\Sigma$ at $x$ and $f_*$ is the
derivative of $f$) extends
continuously across the branch points. This yields a bundle $Tf$ above
$\Sigma$, called the {\it image tangent bundle} and, by taking orthogonal
complements a {\it normal bundle $Nf$}. Note that if $M$
is $4$-dimensional and oriented, then $Nf$ is an oriented $2$-plane bundle.\\
\subsubsection{Degrees of the tangent and normal bundles}
NOTATION. If $L$ is a complex line bundle (or equivalently an oriented real
$2$-plane bundle) on an oriented closed Riemann surface $\Sigma$ without
boundary, we
will denote its degree by $c_1(L)$: that is, we use the same notation for the
first Chern class and for its representative in the second integral cohomology group
of $\Sigma$.\\

We consider here an oriented closed surface without boundary
$\Sigma$ and a map $f$ from $\Sigma$ to an oriented $4$-manifold
$M$. The classical Riemann-Hurwitz formula ([G-H]) generalizes to
\begin{prop}([Gau]).
Suppose $f$ is an immersion with branch points $p_1,...,p_k$ of
respective branching orders $m_1,...,m_k$. Then the degree of the
image tangent bundle $Tf$ is given by
$$c_1(Tf)=\chi(\Sigma)+\sum_{i=1}^k m_i$$
where $\chi(\Sigma)$ denotes the Euler caracteristic of $\Sigma$.
\end{prop}
Hence we derive the degree of $Tf$ by taking the formula for the degree of
the tangent
bundle of an immersed surface (which is $\chi(\Sigma)$) and we add a
correction
term for each branch point, namely its branching order.\\

One of the purposes of this paper will be to achieve a similar
formula for the degree of the normal bundle if $M$ is
$4$-dimensional; that is, to estimate the contribution of a branch
point to the degree of the normal bundle of a branched immersion. To
this effect we now  recall the classical formula for the normal
bundle in the non branched case:
\begin{prop}.
Let $\Sigma$ be a closed oriented Riemann surface without boundary, let $M$ be an oriented smooth $4$-manifold and
let $f:\Sigma\longrightarrow M$ be an immersion with only transverse double
points.
Then the degree of the normal bundle $Nf$ is given by
$$[f(\Sigma)].[f(\Sigma)]-2D_f$$
where $[f(\Sigma)].[f(\Sigma)]$ denotes the self-intersection number of the
$2$-homology class of $f(\Sigma)$ and $D_f$ is the signed number of double
points
of $f$.
\end{prop}
\section{The result}
\begin{thm}
We consider $\Sigma_1$,...$\Sigma_n$ compact connected oriented Riemann surfaces without boundary and $M$ an oriented smooth $4$-manifold. For each
$i=1,...,n$ we let
$$f_i:\Sigma_i\longrightarrow M$$
be a branched immersion and we denote by $Nf_i$ the corresponding normal
bundle. We put ${\mathcal S}=\cup_{i=1}^n \Sigma_i$ and we define a map $f:{\mathcal S}\longrightarrow M$ be imposing its
restriction to each $\Sigma_i$ to be equal to $f_i$.\\

We assume $f({\mathcal S})$  has isolated sigular points  and we let
$[f({\mathcal S)}]$ be its $2$-homology class in $M$.\\

We endow $M$ with a Riemannian metric and denote by $\Bbb{S}_\epsilon(p)$ the sphere centered at a point $p$ in $M$ and of radius $\epsilon$.\\
The link $\Bbb{S}_\epsilon(p)\cap f({\mathcal S})$ is a disjoint union of closed braids
$\Gamma^\epsilon_1,...,\Gamma^\epsilon_s$ which have the same axis up to orientation. Denoting by $e(\Gamma_i^\epsilon)$
the algebraic crossing number of $\Gamma_i^\epsilon$ and by $lk(\Gamma^\epsilon_i,\Gamma_j^\epsilon)$ the linking number of
 $\Gamma^\epsilon_i$ and $\Gamma_j^\epsilon)$ for $i\neq j$, we put
$$E(p)=\sum_{i=1}^s e(\Gamma_i^\epsilon)+ 2\sum_{1\leq i<j\leq s}lk(\Gamma^\epsilon_i,\Gamma_j^\epsilon)$$
Let $p_1,...,p_k$ be the singular points of $f({\mathcal S})$. We
have
$$\sum_{i=1}^n degree(Nf_i)=[f({\mathcal S})].[f({\mathcal S})]-\sum_{m=1}^k E(p_m)$$
where $[f({\mathcal S})].[f({\mathcal S}])$ denotes the self-intersection number of
$f({\mathcal S})$.
\end{thm}
REMARK. Since we mention the linking number of the $\Gamma_i^\epsilon$'s, we need to specify their orientation. We denote
by $\Bbb{B}_\epsilon(p)$ the ball in ${\mathcal S}$ of radius $\epsilon$. Then each $\Gamma_i^\epsilon$ is the image {\it via}
$f$ of $\partial U_i^\epsilon$ where $U_i^\epsilon$ is contained in ${\mathcal S} cap f^{-1}(\Bbb{B}_\epsilon(p))$. The $2$-dimensional
surface $U_i^\epsilon$ inherits the orientation of ${\mathcal S}$: this, in turn, yields an orientation for $\partial U_i^\epsilon$ and
thus for $\Gamma_i^\epsilon$.

\subsection{The case of a single branched disk}
\subsubsection{The knot of the singularity}
We first prove the theorem in the case where there is no more than one branched disc going through each singular point.
We consider a branched disc $f:D\longrightarrow\Bbb{R}^4$ as in Def. 2.
We denote by $e_i, i=1,...,4$ the values of $\frac{\partial}{\partial x_i}$
at the origin. Possibly after replacing $e_4$ by $-e_4$ we assume the basis
$(\frac{\partial}{\partial x_i})$ to be positive w.r.t. the orientation of $\Bbb{R}^4$. We define a  scalar
product $g_0$ on
$\Bbb{R}^4$  by requiring the basis  $(e_1,e_2,e_3,e_4)$ of $\Bbb{R}^4$
to be orthonormal. We also define a complex structure $J_0$ on
$\Bbb{R}^4$ by setting
$$J_0(e_1)=e_2\ ,\ \  J_0(e_3)=e_4.$$

Finally if $\epsilon$ is a small positive number, we let
$\Bbb{S}_\epsilon$ be the sphere (for the norm $g_0$) centered at
$f(0)$ and of radius $\epsilon$. The complex structure $J_0$ yields
isomorphic contact structures on the $\Bbb{S}_\epsilon$'s (we will
denote all these contact structures by $\xi$).
\begin{prop}
Let $f:D\longrightarrow\Bbb{R}^4$ be a branched disc as in Def. 2. Assume moreover that $f$ is a topological embedding in a neighborhood of the origin. There
exists
a number $R$, $0<R<1$ such that for a
small enough positive number $\epsilon$, the curve defined by
$$K^\epsilon=\Bbb{S}_\epsilon\cap f(D(0,R))$$
is a closed braid with braid index $N$. For different $\epsilon$'s the
$K^\epsilon$'s have the same transverse knot type.\\
\end{prop}
It follows that the $K^\epsilon$'s have the same Bennequin self-linking number; since their braid index is the same, they also have the same algebraic
crossing number which we will denote by $e(K)$. This quantity appears in
\begin{lem}
Let $\Sigma$ be an oriented closed surface  without boundary, let $M$ be an
oriented $4$-manifold and let $f:\Sigma\longrightarrow M$ be an immersion which has transverse double points and also branch points $p_1,..., p_k$. \\
Assume that there is only one branched disc going through each $p_i$ and that
the $p_i$'s are isolated singularities of $f(\Sigma)$. \\

For each $i$ we denote by $K_{p_i}$ the closed braid defined by the
branch point $p_i$ as explained above. Then the degree of the normal
bundle $Nf$ is given by
$$[f(\Sigma)].[f(\Sigma)]-2D_f+\sum_{i=1}^{k}e(K_{p_i})$$
(for the notations, see Prop.2).
\end{lem}
Note the formal similarity with Prop. 1.
\subsection{The braids; proof of Prop. 3}
We give a proof of Prop. 3, closely inspired by [Mi]. We denote by $\|.\|$ the norm on $\Bbb{R}^4$ defined by $g_0$. Going back to
Def. 2, we derive
\begin{lem}
There exists a  number $R$, $0<R<1$ such that the following is true: \\
1) The only critical point of
$z\mapsto\|f(z)\|$ in $D(0,R)$ is $0$\\
2) $f(D(0,R)-\{0\})$ is a smooth submanifold of $\Bbb{R}^4$.
\end{lem}
Thus, for $\epsilon>0$ small enough $K^\epsilon$ is smooth. The plane tangent to $K^\epsilon$ at a point $f(z)$ of $K^\epsilon$ is the intersection of the
tangent space to $\Bbb{S}_\epsilon$ at $f(z)$ and the tangent plane to the surface
$f(D)$ at $f(z)$. This last plane, which we  denote $T_{f(z)}f(D)$ is generated by the vectors $\frac{\partial f}{\partial x}, \frac{\partial f}{\partial y}$.
Close to the origin, $f$ looks very much like a holomorphic function, that is,
we derive from  the expression of $f$ in Def. 1,
$$\|\frac{\partial f}{\partial x}\|= |z^{N-1}+o(|z^{N-1}|),\ \
\frac{\partial f}{\partial y}=J_0\frac{\partial f}{\partial x}+o(|z^{N-1|}).$$
It follows that the $K^\epsilon$'s are all transverse knots.\\

The reader can also see that the $K^\epsilon$'s are closed braids and that their axis is the great circle in the plane generated by $e_3, e_4$.\\
To prove that the $K^\epsilon$'s are transversally isotopic we will now recall
a construction from [Mi]. Milnor constructs a vector field $X$ in a
small ball $B$ in $\Bbb{R}^4$ centered at the origin with the
following
properties:\\
1) $X$ is everywhere tangent to $f(D)$\\
2) for every point $p$ in $B$ distinct from the origin, we have
$$<X(p),p>>0.$$
Thus going along the integral curves of $X$ will give an isotopy between
$K^{\epsilon'}$ and $K^{\epsilon}$ for $0<\epsilon'<\epsilon$. \\

We need to point out that Milnor's construction is about {\it
holomorphic functions}, that is, he assumes the map $f$ in Def. 2 to
be holomorphic. However for this specific part of his book (i.e. the
construction of the vector field) holomorphicity is not necessary
and the only thing  needed for the
construction to work is Lemma 1 above. This concludes the proof.\\

We end this section by lemmas which we will use only in \S 3 below.\\
\begin{lem}
For $\epsilon$ small enough, there is a
function
$r_\epsilon$ into $\Bbb{R}^+$
such that $K^\epsilon$ is parametrized by $f(r_\epsilon(t)e^{it})$.
\end{lem}
PROOF. Apply the implicit function theorem to the function
$$(r,\theta)\mapsto <f(re^{i\theta}), f(re^{i\theta})>.$$
We derive
\begin{lem}
We let $T_0f$ be the image tangent plane to $f$ at $0$. Let $Q$ be a
$2$-plane in $\Bbb{R}^4$ which verifies
$$(T_0f)^\perp\cap Q=\{0\}\ \ \ \ (1)$$
(where $P^\perp$ denotes the orthogonal complement of $P$ in $\Bbb{R}^4$). For
$\epsilon>0$, we put $Q^\epsilon=Q\cap \Bbb{S}_\epsilon$. For $\epsilon$ small enough, $Q^\epsilon$ is a braid axis (up
to orientation) for the knot $K^\epsilon$ in $\Bbb{S}_\epsilon$.
\end{lem}
PROOF. We let $\pi_Q:\Bbb{R}^4\longrightarrow Q$ be the orthogonal projection and
we put $\pi_Q(e_1)=\epsilon_1, \pi_Q(e_2)=\epsilon_2$.
Note that $\epsilon_1$ and $\epsilon_2$ are independent vectors. For $z=re^{it}$ in $D$, we have
$$\pi_Q(f(z))=r^N\cos(Nt)\epsilon_1+ r^N\sin(Nt)\epsilon_2+v(z)$$
where $v(z)=o_1(|z|^N)$. Since the $\epsilon_i$'s are independant, there exists a real number $\eta>0$ such that
$$\forall r\in \Bbb{R}, t\in [0,2\pi],
\|r^N\cos(Nt)\epsilon_1+ r^N\sin(Nt)\epsilon_2\|\geq\eta r^N\ \ \ (2).$$
The projection of $\Gamma^\epsilon$ to $Q$ writes as
$$r_\epsilon(t)^N\cos(Nt)\epsilon_1+ r_\epsilon(t)^N\sin(Nt)\epsilon_2+
v(r_\epsilon(t)).$$
We derive Lemma 5 using Lemma 1 2) and (2).
\begin{lem}
The number $e(K^\epsilon)$ does not depend on the braid axis we consider.
\end{lem}
PROOF. We work in the Grassmannian of oriented $2$-planes $G_2^+(\Bbb{R}^4)$.
We recall its canonical splitting into the product of the unit spheres of the eigenspaces of
$\Lambda^2(\Bbb{R}^4)$ (the reader unfamiliar
with these notations can look at \S 2.3.2),
$$G_2^+(\Bbb{R}^4)=\Bbb{S}(\Lambda^+(\Bbb{R}^4))\times \Bbb{S}(\Lambda^-(\Bbb{R}^4)).$$

We fix an orientation of $T_0f$ and write
$T_0f=\frac{1}{\sqrt{2}}(H+K)$, with $H\in
\Bbb{S}(\Lambda^+(\Bbb{R}^4))$ and
$K\in\Bbb{S}(\Lambda^-(\Bbb{R}^4))$. Writing the plane $Q$ as a
$2$-vector and denoting by $<,>$ the scalar product on
$\Lambda^2(\Bbb{R}^4))$ induced by the scalar product on
$\Bbb{R}^4$, (1) translates as
$$<Q,H+K>\neq 0.$$
Hence the set ${\mathcal U}$ of planes verifying (1) is the
complement of a hypersurface in the Grassmannian,
and has two connected components.\\

If we denote by $T_0f^+$ and $T_0f^-$ the tangent plane $Tf$ endowed
with the opposite orientations, we see that they belong to different
components of ${\mathcal U}$. On the other hand, they define the
same crossing number for $K^\epsilon$ (the reader can check this by
looking at a picture in $\Bbb{R}^3$ and figuring out how we get the
same crossing number if we take $\Bbb{R}e_3$ or $\Bbb{R}(-e_3)$ as
a braid axis for a knot.\\
\subsection{Proof of Lemma 2: one branch disk per branch point}
\subsubsection{The idea of the proof}
Since there is a formula for the normal bundle of immersed surfaces with only
transverse double points, we will try and deform $f$ into such an immersion
$\tilde{f}$ while keeping it fixed outside a neighbourhood of the branch points. For technical reasons we do not seek a $\tilde{f}$ which has the same degree
for its normal bundle as $f$. Instead we will require $\tilde{f}$ to verify
$$c_1(Tf)+c_1(Nf)=c_1(T\tilde{f})+c_1(N\tilde{f})\ \ \  (*)$$

Since we know that $c_1(Tf)-c_1(T\tilde{f})$ is equal to the global
branching order of $f$ (the sum of the branching orders of all the
branch points), the
knowledge of the signed number of double points of $\tilde{f}$ will be enough for us to derive the degree of $Nf$.\\

We will explain later how we relate the number of double points of  $\tilde{f}$ to the algebraic crossing number of the $K_{p_i}$'s. For the moment let us just say that the key ingredient is another fact from Bennequin's paper which we now recall.\\
We denote by $\Gamma$ the circle in $\Bbb{S}^3$ defined by the equations
$$x_3=x_4=0.$$
If $K$ is a closed braid $C^1$-close enough to $\Gamma$, and $X$ is a non-zero vector in the plane defined by the equations $x_1=x_2=0$, we denote by $\tilde{K}$ a
knot obtained by pushing $K$ slightly in the direction of $X$. Then the algebraic crossing number $e(K)$ is the linking number $lk(K,\tilde{K})$ between $K$ and $\tilde{K}$.
\subsubsection{Preliminaries: the twistor approach}
1) {\bf A Gauss map}\\
If $f$ is a immersion (possibly with branched points) from a Riemman surface
$\Sigma$ to a $4$-manifold $M$, the quantity $c_1(Tf)+c_1(Nf)$ is the degree of a complex line bundle above $\Sigma$; this bundle can be seen as the
pull-back of another complex line bundle above the {\it twistor space} $Z^+(M)$ of $M$ {\it via} the lift of $f$ into $Z^+(M)$. For this reason we recall now the basic facts about twistor spaces.\\
Consider $\Bbb{R}^4$ endowed with a scalar product. We denote by
$\Lambda^2(\Bbb{R}^4)$ the space of $2$-vectors of $\Bbb{R}^4$ (or equivalently exterior $2$-forms). It admits an involution, the Hodge operator $*$ which
can be defined by requiring:\\
{\it if $(\epsilon_1, \epsilon_2, \epsilon_3, \epsilon_4)$ is a positive orthonormal basis of $\Bbb{R}^4$
then }
$$*(\epsilon_1\wedge \epsilon_2)=\epsilon_3\wedge \epsilon_4.$$
The $+1$-eigenspace of $\Lambda^2(\Bbb{R}^4)$ w.r.t. $*$ is denoted by
$\Lambda^+(\Bbb{R}^4)$. It is a $3$-plane; the metric and orientation on
$\Bbb{R}^4$ yield a metric and orientation on $\Lambda^+(\Bbb{R}^4)$ as well.
We denote by $Z^+(\Bbb{R}^4)$ the unit sphere of $\Lambda^+(\Bbb{R}^4)$; it
inherits from the $SO(3)$-structure on $\Lambda^+(\Bbb{R}^4)$ a natural
$U(1)$-structure.\\

Let $\Sigma$ be a Riemann surface and let
$f:\Sigma\longrightarrow\Bbb{R}^4$ be an immersion with branch
points. We mentioned above the Gauss map from $\Sigma$ to the
Grassmannian of oriented $2$-planes in $\Bbb{R}^4$; we can derive a
map (which we also call a Gauss map):
$$\hat{f}:\Sigma\longrightarrow Z^+(\Bbb{R}^4)$$
$$x\mapsto \frac{1}{\sqrt{2}}(f_*(T_x\Sigma)+*f_*(T_x\Sigma))$$
Note: in this writing  a plane $P$ is viewed as a $2$-vector (i.e. the $2$-vector
$\epsilon_1\wedge\epsilon_2$ where $(\epsilon_1, \epsilon_2)$ is a positive orthonormal basis of $P$).\\
The image tangent (resp. normal) vector bundle $Tf$ (resp. $Nf$) admits a
natural $U(1)$-structure and we have the following isomorphism of complex
line bundles
\begin{prop}
$$Tf\otimes_{\Bbb{C}} Nf\cong \hat{f}^*TZ^+(\Bbb{R}^4).$$
\end{prop}

2) {\bf Twistor bundles}\\
If we replace $\Bbb{R}^4$ by a Riemannian $4$-manifold $M$ these constructions become
bundle constructions. We define the bundles $\Lambda^2(M)$, $\Lambda^+(M)$
and $Z^+(M)$. $Z^+(M)$ is a bundle of oriented $2$-spheres above $M$. As such
it admits a {\it vertical tangent bundle} $T^v$; note that $T^v$ is a
$U(1)$-bundle.\\

If $f:\Sigma\longrightarrow M$ is an immersion with branch points,
we can define its {\it twistor lift} $\hat{f}:\Sigma\longrightarrow
Z^+(M)$ similarly as above. The metric on $M$ together with the
orientations on $\Sigma$ and $M$ yield $U(1)$-structures on the
image tangent (resp. normal) bundle $Tf$ (resp. $Nf$) and  Prop. 4
becomes
\begin{prop}
There is an isomorphism of complex line bundles
$$Tf\otimes_{\Bbb{C}}Nf\cong \hat{f}^*T^v.$$
\end{prop}

REMARK. The reader can find a proof of Prop. 4 and 5 in [E-S] and [Vi] with
the following {\it caveat}: the complex structure these two papers assume on
$Z^+(\Bbb{R}^4)$ is conjugate to the one we consider here. In the present paper
we work with the so called {\it natural complex structure} on $Z^+(\Bbb{R}^4)$. Eells and Salamon, however, had to consider the conjugate of the natural
complex structure  in order to get holomorphicity of the twistor lifts of
conformal harmonic maps.

\subsubsection{A family of immersions}
It is enough to prove Lemma 2 in the case where $f$ has a single branch
point $p$. We let $z$ be a holomorphic coordinate in $\Bbb{C}$ around $0$ and we let $(x_1,...,x_4)$ be coordinates in a neighbourhood $U$ of $f(0)=p$ in $M$. In this neighbourhood, the coordinates of $f(z)$ write
$$x_1=Re(z^N)+u_1(z)$$
$$x_2=Im(z^N)+u_2(z)$$
$$x_3=h_1(z)$$
$$x_4=h_2(z)$$
where $u_1(z)$, $u_2(z)$, $h_1(z)$ and $h_2(z)$ are $o_1(|z^N|)$'s.\\

We set out to construct a family of immersions which approach $f$,
coincide with $f$ outside of a neighbourhood of the branch point and
verify condition (*)
above.\\

For a positive real number $r$, we define a cut-off function
$\zeta_r:[0,1]\mapsto [0,1]$ such that
$$\zeta_r(t)=1\ \text{if}\  t<\frac{r}{2}$$
$$\zeta_r(t)=0\ \text{if}\  t>\frac{2}{3}r.$$
For a small positive number $r$ and a small number $\lambda$ we define
$$f_{\lambda,r}:D\longrightarrow M$$
$$
z\mapsto \left|
\begin{array}{c}
Re(z^N)+u_1(z)\\
Im(z^N)+u_2(z)\\
h_1(z)+\lambda\zeta_r(|z|)Re(z)\\
h_2(z)+\lambda\zeta_r(|z|)Im(z)
\end{array}
\right.
$$
Possibly by restricting the disc $D$, we can assume that the function
from $D$ to $\Bbb{R}^2$ defined by
$$z\mapsto (Re(z^N)+u_1(z),Im(z^N)+u_2(z))$$
is an immersion except at $0$. Thus $f_{\lambda,r}$ is an immersion
everywhere. Note that the tangent plane $f_{\lambda,r}^*(T_0 D)$
coincides with the
{\it normal} plane $N_0 f$.\\
$f_{\lambda,r}$ coincides with $f$ in a neighbourhood of the boundary of $D$. By setting $f_{\lambda,r}$ to coincide with $f$ outside of $D$ we extend it to
an immersion from the entire Riemann surface $\Sigma$ to $M$. The surfaces
$f(\Sigma)$ and $f_{\lambda,r}(\Sigma)$ are evidently cohomologous in $M$.
\begin{lem}
There exist a real number $r_0$ and a function $\lambda(r)$ taking positive
values such that: \\
if  $r<r_0$ and $\lambda<\lambda(r)$ , we have
$$c_1(Tf_{\lambda,r})+c_1(Nf_{\lambda,r})=c_1(Tf)+c_1(Nf) $$
\end{lem}
PROOF. We choose a Riemannian metric $g$ on $M$ which satisfies the
requirement:\\
in the above mentioned neighbourhood $U$ of $f(p)$,  the basis of tangent
vectors
$(\frac{\partial}{\partial x_1},...,\frac{\partial}{\partial x_4})$ is
orthonormal w.r.t. $g$. We
denote by $\Lambda^+(M)$ (resp. $Z^+(M)$) the
bundle of $2$-vectors which are $+1$-eigenvectors for the Hodge operator (resp. the twistor space) on $(M,g)$. \\
In view of Prop.5, Lemma 2 will follow from
\begin{slem}
Let $\hat{f}$ (resp. $\hat{f}_{\lambda,r}$) be the lift of $f$ (resp.
$f_{\lambda,r}$) in $Z^+(M)$. Then there exists a real number $r_0$ and a
function $\lambda(r)$ such that:\\
if $r<r_0$ and $\lambda<\lambda(r)$, then $\hat{f}(\Sigma)$ and
$\hat{f}_{\lambda,r}(\Sigma)$ have the same homology class in $Z^+(M)$.
\end{slem}
PROOF OF SUBLEMMA 1. Since these maps coincide outside a neighbourhood of the branch point, it is
enough to focus our attention near the branch point. We need to introduce some
more notations.\\

We can trivialize $\Lambda^+(M)$ above $U$ with the help of the
following orthonormal basis $(H_1,H_2,H_3)$ (the fibres of
$\Lambda^+(M)$ inherit a scalar product from the metric $g$ on $M$)
:
$$H_1=\frac{1}{\sqrt{2}}(e_1\wedge e_2+ e_3\wedge e_4)$$
$$H_2=\frac{1}{\sqrt{2}}(e_1\wedge e_3+e_4\wedge e_2)$$
$$H_3=\frac{1}{\sqrt{2}}(e_1\wedge e_4+e_2\wedge e_3).$$
where $e_i=\frac{\partial}{\partial x_i}$ for $i=1,...,4$.\\
Above $U$, $Z^+(M)$ will be trivialized by the unit sphere of
the space generated by $(H_1,H_2,H_3)$. \\

We now define, for a function $h:D\longrightarrow M$, with $h(0)=p$
the following lift of $h$ in $\Lambda^+(M)$ in $U$
$$\Lambda^+(h)=\frac{\partial h}{\partial x}\wedge
\frac{\partial h}{\partial y}
+*(\frac{\partial h}{\partial x}\wedge \frac{\partial h}{\partial y})$$
where $x$ (resp. $y$) denote the real (resp. complex) part of the coordinate
$z$ on $D$ and $*$ is the Hodge operator on $\Lambda^2(M)$ associated with
the metric $g$.\\

If $h$ is an immersion, its lift in $Z^+(M)$ writes
$$\hat{h}=\frac{\Lambda^+(h)}{\|\Lambda^+(h)\|}.$$
We can now return to the proof of Sublemma 1. Using the above notations, we see that
$$\hat{f}(0)=H_1.$$
Hence there exists a number  $r_1$ such that \\
 $$if\ \ \  |z|<r_1,
<\hat{f}(z),H_1>>0$$
where $<,>$ denotes the scalar product in the fibres of $\Lambda^+(M)$. Thus
Sublemma 1 will follow from
\begin{slem}
There exist a positive real number $r_0$ and a function $\lambda(r)$ such that if $r<r_0$ and $\lambda<\lambda(r)$,\\
$$\forall z\in D(0,r)
 <\Lambda^+(f_{\lambda,r})(z),H_1>>0.$$
\end{slem}
PROOF OF SUBLEMMA 2. We will use the following formula
\begin{slem}
Let $H:D\longrightarrow M$, with $H(0)=p$ which writes in the coordinate system
$(x_i)$ \\
$$h(z)=(Re(z^N)+u_1(z),Im(z^N)+u_2(z),v_1(z),v_2(z))$$
where $u_1, u_2$ are
$o_1(|z|^N)$. Then in a neighbourhood of $0$,
$$<\Lambda^+(H)(z),H_1>=N^2|z|^{2(N-1)}+
(\frac{\partial v_1}{\partial x}\frac{\partial v_2}{\partial y}
-\frac{\partial v_1}{\partial y}\frac{\partial v_2}{\partial x})+o(|z|^{2(N-1)})$$
\end{slem}
The proof of Sublemma 3 is an immediate computation which we omit. We
move on to prove Sublemma 2.\\

We treat separately the cases where $|z|<\frac{r}{2}$ and
$|z|>\frac{r}{2}$. A straightforward computation yields
\begin{slem}
There exist two real functions $F_1$ and $F_2$ of one complex variable,
verifying
$$F_1(z)=o(|z|^{N-1}), F_2(z)=o(|z|^{2(N-1)}),$$
such that, for $|z|<\frac{r}{2}$,
$$<(\Lambda^+f_{\lambda,r}(z),H_1>=N^2|z|^{2(N-1)}+\lambda^2+
\lambda F_1(z)+F_2(z).$$
\end{slem}

So there exists an $r_2>0$ such that for $r<r_2$ and $|z|<\frac{r}{2}$,
$$<(\Lambda^+f_{\lambda,r}(z), H_1>>\frac{3}{2}|z|^{2(N-1)}
+\lambda^2-2\lambda|z|^{N-1}$$
$$=(|z|^{N-1}-\lambda)^2+\frac{1}{2}|z|^{2(N-1)}>0.$$

\begin{slem}
There exist three real functions $F_3$, $F_4$ and $F_5$ of one complex
variable, verifying
$$F_3(z)=o(|z|^{2(N-1)}, F_4(z)=o(|z|^{N-1}), F_5(z)=o(|z|^N)$$
such that, for $\frac{r}{2}<|z|<\frac{2}{3}r$
$$<\Lambda^+{f}_{\lambda,r}(z),H_1>=N^2|z|^{2(N-1)}+F_3(z)+\lambda\zeta_r(|z|)F_4(z)$$
$$+\lambda\zeta'_r(|z|)F_5(z)
+\lambda^2\zeta_r(|z|)^2+\lambda^2|z|\zeta_r(|z|)\zeta_r'(|z|).$$
\end{slem}

We leave the proof (straightforward computation again) to the reader.\\

We put $A(r)=\sup\|\zeta_r'\|$. Sublemma 5 yields a $r_3>0$ such
that for $r<r_3$, $\frac{r}{2}<|z|<\frac{2}{3}r$ and $\lambda<1$,
$$<\Lambda^+{f}_{\lambda,r}(z),H_1>>(N^2-1)(\frac{r}{2})^{2(N-1)}
-\lambda(\frac{2}{3}r)^{N-1}$$
$$-\lambda A(r)(\frac{2}{3}r)^N-\lambda^2(\frac{2}{3}r)A(r)$$
$$>(N^2-1)\frac{r}{2}^{2(N-1)}-\lambda[1+2A(r)].$$

It follows from this discussion that Sublemma 2 will be satisfied if
we take
$$r_0=min(r_1,r_2,r_3)$$
$$\lambda(r)=min(1,\frac{(N^2-1)(\frac{r}{2})^{2(N-1)}}{2[2A(r)+1]}).$$
This concludes the proof of Sublemma 1.\\

We fix a $\lambda,r$ satisfying the assumptions of Sublemma 1.
\begin{slem}
There is a neighbourhood $V_1$ of $0$ in $D(0,1)$ and a neighbourhood $V_2$ of the boundary of $D(0,1)$ such that the restriction of $f_{\lambda,r}$ to
$V_1\cup V_2$ is a smooth embedding.
\end{slem}
PROOF. By looking at the expression of $f_{\lambda,r}$, it is clear that there
exists a neighbourhood $V_1$ (resp. $V_2$) of the origin (resp. the boundary
of $D(0,1)$) such that:\\
1) the restriction of $f_{\lambda,r}$ to each of the $V_i$'s,
$i=1,2$, is an
embedding\\
2) if we denote by $f_{\lambda,r}^{1,2}$ the projection of
$f_{\lambda,r}$ to the plane generated by the first two components,
namely $e_1$ and $e_2$, there exists a number $C$ such that
$$\forall z_1\in V_1, |f_{\lambda,r}^{1,2}(z_1)|>C,\ \
\forall z_2\in V_2, |f_{\lambda,r}^{1,2}(z_1)|<C.$$
By transversal approximation we
construct for each couple $(\lambda,r)$ a sequence
$(f_{\lambda,r}^{(n)})$ of immersions from $\Sigma$ to $M$ which converges
$C^1$
to $f_{\lambda,r}$ and such that for every $n$,\\
i) $f_{\lambda,r}^{(n)}$ has only transverse double points\\
ii) $f_{\lambda,r}^{(n)}$ coincides with $f_{\lambda,r}$ on $V_1\cup V_2$.\\
 Since the $f_{\lambda,r}^{(n)}$'s
converge to $f_{\lambda,r}$, it follows that for $n$ large enough
$$c_1(Nf_{\lambda,r}^{(n)})=c_1(Nf_{\lambda,r})$$
$$=[f(\Sigma)].[f(\Sigma)]-2D_f-2D_{\lambda,r}^{(n)}\ \ \ (**)$$
where $D_{\lambda,r}^{(n)}$ denotes the number of double points (counted w.r.t.
sign) of the restriction of $f_{\lambda,r}^{(n)}$ to $D(0,r)$.\\
Since $c_1(Tf)+c_1(Nf)=c_1(Tf_{\lambda,r})+c_1(Nf_{\lambda,r})$, we derive
$$c_1(Nf)=[f(\Sigma)].[f(\Sigma)]-2D_f-(N-1)-2D_{\lambda,r}^{(n)}.$$
In order to proceed and interpret the last two terms of the previous line,
we need to recall a few facts about topology of $4$-manifolds
\subsubsection{Framing of a knot in $\Bbb{S}^3$}
A {\it framing} of a knot $K$ in $\Bbb{S}^3$ is a vector field $Y$ along $K$ on $\Bbb{S}^3$  which is never tangent to $K$.\\
 Since the normal bundle to $K$ in $\Bbb{S}^3$ is trivial, the homotopy classes of framings are parametrized by  $\pi_1(\Bbb{R}^2-\{0\})$, i.e. $\Bbb{Z}$.
We denote by $\tilde{K}$ the knot obtained by pushing $K$ in the direction of
$Y$. We let $S$ (resp. $\tilde{S}$) be a smooth surface {\it embedded} in
$\Bbb{B}^4$
bounded by $K$ (resp. $\tilde{K}$). The number of intersection points of $S$
and $\tilde{S}$ (counted w.r.t. sign) is equal to the linking number of $K$ and
$\tilde{K}$. It is also called the {\it self-intersection number of $S$ w.r.t.
the framing $Y$}. Now suppose that $Y$ is orthogonal to the restriction to $K$ of the tangent bundle to $S$  and that it extends to a section (which we also
denote $Y$) of the normal
bundle $NS$ above $S$. Then the number of zeroes $N(Y)$ of $Y$ on $NS$ is
also equal to $lk(K,\tilde{K})$. Finally if $S$ is not embedded but immersed
with only transverse double points, and $Y$ extends to a section of $NS$ which does not vanish at any double point, we have
$$N(Y)=lk(K,\tilde{K})-2D_S \ \ \ (3)$$
where $D_S$ denotes the number of double points of $S$ counted w.r.t. sign.
\subsubsection{End of the proof}
Since $f(D)$ is contained in $U$ which is isometric to an Euclidean $4$-ball, we can define the projection $X_N$ of the vector $X$ on the normal bundle $Nf$.
Because $X$ and $X_N$ coincide at the point $p$, it follows that for $\epsilon$ small enough, $X_N$ defines a framing of $K^\epsilon$ which is homotopic to $X$
We choose an $r<r_0$ such that $f(D(0,r))$ is contained.in the ball centered at
$f(p)$ and of radius $\epsilon$ and a $\lambda<\lambda(r)$. We fix these values of $r$ and $\lambda$ and we choose an $n$ such that (**) is verified. We denote by
$X^{(n)}_N$ the projection of $X$ to the normal bundle $Nf^{(n)}_{\lambda,r}$;
it coincides with $X_N$ on the boundary $K^\epsilon$. We let $N(X_N^{(n)})$ be
the number of zeroes counted w.r.t. sign . Putting together (3) and the
property of the algebraic crossing number described in \S 2.3.1, we get that
$$N(X_N^{(n)})=e(K^\epsilon)-2D_{f_{\lambda,r}^{(n)}}.$$

In view of (**), we will be able to conclude the proof of Lemma 1
once we have proven the
\begin{lem}
i) For $n$ large enough, $X_N^{(n)}$ only vanishes at $0$\\
ii) The index of $X_N^{(n)}$ at $0$ is equal to $N-1$.
\end{lem}
PROOF. We denote by $X_N(\lambda,r)$ the orthogonal projection of $X$ to the
normal bundle $Nf_{\lambda,r}$. Going back to the explicit formula for
$f_{\lambda,r}$ we notice that $X$ is only tangent to $f_{\lambda,r}(D)$ at $0$. Hence, there exist a positive number $ \alpha$ such that:\\
{\it outside of $D(0,\frac{r}{2})$, the norm $\|X_N(\lambda,r)\|$  is bounded below by $\alpha$}.\\

It follows that for $n$ large enough, $X_N^{(n)}$ does not vanish
outside of $D(0,\frac{r}{2})$. On the other hand, inside
$D(0,\frac{r}{2})$, $X_N^{(n)}$ coincides with $X_N(\lambda,r)$; and
this latter vector field only vanishes at $0$. This
concludes the proof of i).\\

We prove ii) by deriving an explicit formula for $X_N(\lambda,r)$.
In order
to do this, we introduce some notations.\\

We let $P$ be a $2$-plane in $\Bbb{R}^4$ and we write
$P=\epsilon_1\wedge \epsilon_2$ for an orthonormal basis $(e_1,e_2)$
of $P$. We denote by $J_P$ the complex structure on $P$ compatible
with the metric and orientation, that is
$$J_P(\epsilon_1)=\epsilon_2,\ \  J_P(\epsilon_3)=\epsilon_4.$$
We let $*:\Lambda^3(\Bbb{R}^4)\longrightarrow\Lambda^1(\Bbb{R}^4)$ be the Hodge operator, characterized by the following
property: if $(u_1,u_2,u_3,u_4)$ is a
positive orthonormal basis of $\Bbb{R}^4$, then
$$*(u_1\wedge u_2\wedge u_3)=u_4.$$
Finally we denote by $\pi_P$ the orthogonal projection from $\Bbb{R}^4$ to $P$.
It is straightforward to check that, under the above notations,
$$\forall X\in\Bbb{R}^4,\ \ \  \pi_P(X)=-J_P[*(X\wedge *P)].$$
We can trivialize $Nf_{\lambda,r}$ in a neighbourhood of $0$ by orthogonal projection to $Tf(0)$. We notice that
$$e_1=*(e_3\wedge e_2\wedge e_4), e_2=-*(e_3\wedge e_1\wedge e_4)\ \ \ (4)$$

A quick computation yields, for $z=\rho e^{i\theta}$,
$$(\frac{\partial f_{\lambda,r}}{\partial x}\wedge
\frac{\partial f_{\lambda,r}}{\partial y})(z)
=\lambda N\rho^{N-1}
[(\cos(N-1)\theta)e_1\wedge e_4+
\sin(N-1)\theta)e_2\wedge e_4]$$
$$+A+B$$
where $\|B\|$ is a $o(\rho^{N-1})$ and $A$ is a $2$-vector orthogonal to
$e_1\wedge e_2$ and \\$e_2\wedge e_4$. It follows from the formula above for $e_1$ and $e_2$ that $X_N$ writes in the trivialization of $Nf_{\lambda,r}$
around $0$ by $Tf(0)$
$$\|\frac{\partial f_{\lambda,r}}{\partial x}\wedge
\frac{\partial f_{\lambda,r}}{\partial y}\|X_N=
\lambda N\rho ^{N-1}
[(\cos(N-1)\theta)e_1+
\sin(N-1)\theta)e_2]+o(\rho^{N-1}).$$
This concludes the proof of Lemma 1
\section{The general case}
To introduce the case of several disks meeting at the same point, we point out
the simplest example, i.e. two smoothly embedded disks meeting
transversally: the resulting link is the Hopf link $K_1\cup K_2$. We have $e(K_1)=e(K_2)=0$ and $lk(K_1,K_2)=\pm 2$:
we already know that a transverse
double point contributes $\pm 2$ to the degree of the normal bundle. \\

We will now prove Th. 1; as in the previous section, we can assume
that there is  a single singular point $p$.
The point $p$ may have several preimages {\it via} $f$ and even several preimages {\it via} the same $f_k$. \\

For a $k\in\{1,...,n\}$, we denote by $q_{1_k},..., q_{s_k}$ the
preimages of $p$ in $\Sigma_k$ (of course there can be some
$k$'s for which $p$ has no preimage in $\Sigma_k$).\\

For every $i_k$, we introduce a neighborhood $D_{i_k}$ of $q_{i_k}$
in $\Sigma_k$; we take the $D_{i_k}$'s to be all disjoint. On each
$D_{i_k}$ we assume $f_k$ to be parametrized as in Def. 1 (with
$f_k(0)=p$). We denote by $N_{i_k}$ the integer which appears in
that definition. Please note: here, unlike in Def. 1,  we do not
require $N_{i_k}$ to be strictly greater than $1$;
if $q_{i_k}$ is a regular point of $f_k$, then $N_{i_k}=1$.\\

We now consider the sphere $\Bbb{S}_\epsilon(p)$ of radius
$\epsilon$ centered at $p$ and for each $i_k$ we put
$$\Gamma^\epsilon_{i_k}=\Bbb{S}_\epsilon(p)\cap f_\epsilon(D_{i_k}).$$
We put $L^\epsilon=\cup_{i_k}\Gamma^\epsilon_{i_k}$.\\

The results from the previous paragraph show us that for $\epsilon$ small enough\\
1) the isotopy type of  $\Gamma^\epsilon_{i_k}$ does not depend on $\epsilon$\\
2) there exists a $2$-plane $Q$ such that -up to orientation - all the $\Gamma^\epsilon_{i_k}$ are
closed braids of braid axis $Q^\epsilon=Q\cap \Bbb{S}_\epsilon(p)$.\\

The quantity
$$E(L^\epsilon)=\sum_{i_k}e(\Gamma^\epsilon_{i_k})
+2\sum_{i_k, j_{l}}lk( \Gamma^\epsilon_{i_k},\Gamma^\epsilon_{j_l})\ \ \  (5)$$
 does not depend on $\epsilon$ small enough (see 1) above).\\
We now use the constructions of the previous paragraph to introduce for every $i,k$ an immersion
$h_{i_k}:D_{i_k}\longrightarrow M$ with the following properties\\
1) the lift of $h_{i_k}$ into the twistor space $Z^+(M)$ is close enough to the
lift of $f_k$\\
2) $h_{i_k}$ coincides with $f_k$ near the boundary of $D_{i_k}$\\
3) $h_{i_k}$ has only transverse double points. Their number is equal to
$$\frac{1}{2}[e(\Gamma^\epsilon_{i_k})-(N_{i_k}-1)]$$
By transverse approximation we can also assume \\
4) if $\{i,k\}\neq \{j,l\}$, then $h_{i_k}(D_{i_k})$ and
$h_{j_l}(D_{j_{l}})$ meet transversally: the number of their intersection points is
$$lk(\Gamma^\epsilon_{i_k}\Gamma^\epsilon_{j_k}).$$
Putting together the $h_{i_k}$'s, we can construct an immersion
$$h_k:\Sigma_k\longrightarrow M$$
which coincides with $h_{i_k}$ on the $D_{i_k}$, with $f_i$ outside and verifies
$$c_1(Th_k)+c_1(Nh_k)=c_1(Tf_k)+c_1(Nf_k).$$
It follows that
$$c_1(Nf_k)=[f_k(\Sigma_k)].[f_k(\Sigma_k)]-\sum_i e(\Gamma^\epsilon_{i_k})
-2\sum_{i,j}lk(\Gamma^\epsilon_{i_k},\Gamma^\epsilon_{j_k})\ \ \ (6)$$
If $k$ and $l$ are different, $h_k(\Sigma_k)$ and $h_l(\Sigma_l)$ meet
transversally and we have
$$[f_k(\Sigma_k)].[f_l(\Sigma_l)]=[h_k(\Sigma_k)].[h_l(\Sigma_l)]=
2\sum_{i,k,j,l}lk(\Gamma^\epsilon_{i_k},\Gamma^\epsilon_{j_l})$$

On the other hand,
$$[f({\mathcal S})].[f({\mathcal S})]=[\sum_k f_k(\Sigma_k)].
[\sum_k f_k(\Sigma_k)]$$
$$=\sum_k [f_k(\Sigma_k)].[f_k(\Sigma_k)]-
2\sum_{k,l, k\neq l}[f_k(\Sigma_k)].[f_l(\Sigma_l)]$$
$$=\sum_k [h_k(\Sigma_k)].[h_k(\Sigma_k)]-
2\sum_{k,l,k\neq l}[h_k(\Sigma_k)].[h_l(\Sigma_l)]$$
$$=\sum_k [h_k(\Sigma_k)].[h_k(\Sigma_k)]-
2\sum_{k,l,k\neq l}\sum_{i_k,j_l}lk(\Gamma^\epsilon_{i_k}, \Gamma^\epsilon_{j_l}).\ \ \ (7)$$
Putting together (5), (6) and (7) yields the result.
\section{Appendix: the self-linking number of iterated torus knots}
In this section we look at a case where the computations of the previous section can be carried explicitely. We consider a branched disk of the following type
$$
z\mapsto\left|
\begin{array}{c}
x_1+ix_2=z^N\\
x_3+ix_4=\sum^s_{j=0}(a_j z^{\mu_j}+b_j\bar{z}^{\mu_j})=P(z)
\end{array}\right. \ (***).
$$
with the condition
$$\forall i,\ \ \  a_i b_i=0.$$

The integers $N,s,\mu_0,...,\mu_s$ verify
$$0<N<\mu_0<...<\mu_s.$$

We assume moreover that $gcd(N,\mu_0,...,\mu_s)=1.$ Our attention
was drawn to this type of branch point by Micallef and White's work
([M-W]) on area minimizing surfaces (i.e. surfaces which are a local
minimum of the area for every small deformation). They proved that
the
branch point of such maps are of the type (***) above.\\

If $\epsilon$ is a small enough real positive number we denote by
$\Bbb{S}_\epsilon$ the sphere centered at the origin and of radius
$\epsilon$ and we put $K^\epsilon=f(D)\cap \Bbb{S}_\epsilon$. For
$\epsilon$ small enough, $K^\epsilon$ is an iterated torus knot. It
is well-known how to compute the self-linking number of these knots,
for example see the formula given in [B-W]. However we would like to
show here that this number can also be computed using the branched
immersion approach:
\begin{prop}
Using the notations above and for $\epsilon$ small enough, we have
$$sl(K^\epsilon)=\sum_{j=0}^s(Q_j-Q_{j+1})\tau_j$$
where the $Q_{j+1}$'s and $\tau_j$'s are defined by\\
1) $N=Q_0$\\
2) $Q_{j+1}$ is the greatest common divisor of $Q_j$ and $\mu_j$\\
3)$$\tau_j=\left\{
\begin{array}{c}
\mu_j-1\ \ if\ \ a_j\neq 0\\
-(\mu_j+1)\ \ if\ \ b_j\neq 0
\end{array}\right.
$$
\end{prop}
REMARK. Note the formal similarity with the formula for the Milnor number ([Mi], p.93).\\

Since this result is already known we just outline the main steps of
the proof. We introduce a cut-off function $\zeta$ and define maps
$f_{\lambda,r}$'s as in \S 2.3.3 above. It will turn out that they
only have transverse double points so the (signed) number of their
double points will yield the self-linking number of $K^\epsilon$. To
derive the double points we put
 $S_\nu=P(z)-P(\nu z)+\lambda z-\lambda\nu z$: a double point is a
pair $\{z,\nu z)\}$ with $S_\nu(n)=0$. So our goal is to compute
the cardinal of
$${\mathcal D}=\{z\in\Bbb{C}\ \  |z|<\frac{r}{2}\ \ \  S_\nu(z)=0\}.$$
We might get the same double point {\it via} two different $\nu$'s, so we state
\begin{lem}
1) There is a $1-1$ correspondance
$${\mathcal D}_\nu\longrightarrow{\mathcal D}_{\nu^{-1}}$$
$$z\mapsto\nu z.$$
2) If $\{z,\nu z\}=\{z',\nu z'\}$ is a double point, there either $z=z'$,
$\nu=\nu'$ or $z'=\nu z, \nu\nu'=1$.
\end{lem}

We next compute the cardinal of ${\mathcal D}_\nu$ for a given $\nu$
(without paying attention to their sign). We separate the $N$-th
roots of $1$ into different classes, namely we set for $j=0,...,s$
$${\mathcal R}_j=\{\nu\in\Bbb{C},\ \  \nu^{Q_j}=1\ \ \ and \ \
\nu^{Q_{j+1}}\neq 1\}$$

If $\nu$ belongs to a ${\mathcal R}_j$, then for every $i\leq j$, we
have
$$\mu^i=1.$$

It follows that, for $\nu\in {\mathcal R}_j$, $S_{nu}(z)$ is of one
of the
following two forms depending on which of $a_{j+1}$ and $b_{j+1}$ is non zero:\\
1) if $a_{j+1}\neq 0$, $S_\nu(z)=a_{j+1}(1-\nu^{\mu_{j+1}})z^{\mu_j}
+\lambda(1-\nu)z+h(z)$\\
2)if $b_{j+1}\neq 0$, $S_\nu(z)=b_{j+1}(1-\bar{\nu}^{\mu_{j+1}})\bar{z}^{\mu_j}
+\lambda(1-\nu)z+h(z)$\\
where in both cases, $h(z)$ is a polynomial in $z$, $\bar{z}$ with all its terms of total degree higher than $\mu_{j+1}$. It is easy to derive that \\
{\it $card({\mathcal D}_\nu)=\mu_j-1$ (resp. $card({\mathcal D}_\nu)=\mu_j+1$}.\\

This being established, the  last part in the proof consists in
checking that the sign of a double point in ${\mathcal R}_j$ is
positive (resp. negative) if $a_{j+1}\neq 0$ (resp. $b_{j+1}\neq
0$).

Northeastern University, Department of Mathematics, \\360 Huntington
Avenue, Boston MA 02115, USA \\
mville@math.jussieu.fr

\end{document}